\documentclass[12pt,a4paper,reqno]{amsart}  
\usepackage{amssymb} \usepackage{amscd} \usepackage{times}
\newcommand{\be} {\begin{equation}}
\newcommand{\ee} {\end{equation}}
\newcommand{\bea} {\begin{eqnarray}}
\newcommand{\eea} {\end{eqnarray}}
\newcommand{\Bea} {\begin{eqnarray*}}
\newcommand{\Eea} {\end{eqnarray*}}

\def\zbb{\mathbb{Z}}  
  
  \def\phi{\varphi}
 \def\p1{{\mathbb{P}^1_\zbb}}

\newtheorem{Theorem}{\quad Theorem}[section]

\newtheorem{Corollary}[Theorem]{\quad Corollary}
\newtheorem{Lemma}[Theorem]{\quad Lemma}
\newtheorem{Proposition}[Theorem]{\quad Proposition}

\begin{document}
\title{ An estimate on Riemannian manifolds of dimension 4.}
\author{Samy Skander Bahoura}
\address{Departement de Mathematiques, Universite Pierre et Marie Curie, 2 place Jussieu, 75005, Paris, France.}
\email{samybahoura@yahoo.fr}
\date{\today}
\maketitle
\begin{abstract}
We give an estimate of type $ \sup \times \inf $ on Riemannian manifold of dimension 4 for a Yamabe type equation.
\end{abstract}
\section{Introduction and Main Results}
In this paper, we deal with the following Yamabe type equation in dimension $ n= 4 $:
 \be \Delta_g u+h u=8 u^{3},\,\, u >0, \,\,   \label{(1)}\ee

Here, $ \Delta_g $ is the Laplace-Beltrami operator and $ h $ is an arbitrary bounded function..

The equation $ (\ref{(1)}) $ was studied a lot, when $ M =\Omega \subset {\mathbb R}^n $ or $ M={\mathbb S}_n $ see for example, [2-4], [11], [15]. In this case we have a $ \sup \times \inf $ inequality.
The corresponding equation in two dimensions on open set $ \Omega $ of $ {\mathbb R}^2 $, is:
 \be \Delta u=V(x)e^u, \label{(2)} \ee
The equation $ (\ref{(2)}) $ was studied by many authors and we can find
very important result about a priori estimates in [8], [9], [12],
[16], and [19]. In particular in [9] we have the following interior
estimate:
 $$  \sup_K u  \leq c=c(\inf_{\Omega} V, ||V||_{L^{\infty}(\Omega)}, \inf_{\Omega} u, K, \Omega).  $$
And, precisely, in [8], [12], [16], and [20], we have:
 $$ C \sup_K u + \inf_{\Omega} u \leq c=c(\inf_{\Omega} V, ||V||_{L^{\infty}(\Omega)}, K, \Omega),  $$
and,
 $$ \sup_K u + \inf_{\Omega} u \leq c=c(\inf_{\Omega} V, ||V||_{C^{\alpha}(\Omega)}, K, \Omega).  $$
where $ K $ is a compact subset of $ \Omega $, $ C  $ is a positive constant which depends on $\dfrac{\inf_{\Omega} V}{\sup_{\Omega} V} $, and,  $ \alpha \in (0, 1] $.
When $ 6h=R_g $ the scalar curvature, and $ M $ compact, the equation $(\ref{(1)}) $ is Yamabe equation. T. Aubin and R. Schoen have proved the existence of solution in this case, see for example [1] and [14] for a complete and detailed summary.
When $ M $ is a compact Riemannian manifold, there exist some compactness result for equation  $ (\ref{(1)}) $ see [18]. Li and Zhu see [18], proved that the energy is bounded and if we suppose $ M $ not diffeormorfic to the three sphere, the solutions are uniformly bounded. To have this result they use the positive mass theorem.
Now, if we suppose $ M $ Riemannian manifold (not necessarily compact) 
 Li and Zhang [17] proved that the product $ \sup \times \inf $ is bounded. Here we extend the result of [5].
 Our proof is an extension Li-Zhang result in dimension 3, see [3] and [17], and,  the moving-plane method is used to have this estimate. We refer to Gidas-Ni-Nirenberg for the  moving-plane method, see  [13]. Also, we can see in [3, 6, 11, 16, 17, 10], some applications of this method, for example an uniqueness result. We refer to [7] for the uniqueness result on the sphere and in dimension 3.
Here, we give an equality of type $ \sup \times \inf $ for the equation $ (\ref{(1)}) $ in dimension 4.
In dimension greater than 3 we have other type of estimates by using moving-plane method, see for example [3, 5].
There are other estimates of type $ \sup + \inf $ on complex Monge-Ampere equation on compact  manifolds, see [20-21] . They consider, on compact Kahler manifold $ (M, g) $, the following equation:
 \be 
 \begin{cases}
(\omega_g+\partial \bar \partial \phi)^n=e^{f-t\phi} \omega_g^n, \\
\omega_g+\partial \bar \partial \phi >0 \,\, {\rm on } \,\,  M  \\
\end{cases}  
  \label{(3)} \ee

And, they prove some estimates of type $ \sup_M+m \inf_M \leq C $ or $  \sup_M + m \inf_M \geq C $ under the positivity of the first Chern class of M.
Here, we have,

\begin{Theorem} For all compact set $ K $ of $ M $, there is a positive constant c, which depends only on, $ h_0=||h||_{L^{\infty}(M)}, K, M, g $ such that:
$$ (\sup_K  u)^{1/3} \times \inf_M u \leq c,  $$
for all $ u $ solution of $ (\ref{(1)}) $.
\end{Theorem}
This theorem extend to the dimension 4 a result of the author and of Li and Zhang result, see [17] . Here, we use a different method than the method of Li and Zhang in [17] . Also, we extend a result of [5].

\begin{Corollary}  For all compact set $ K $ of $ M $ there is a positive constant c, such that:

 $$ \sup_K  u \leq c=c(g, m, h_0, K, M) \,\, {\rm if } \,\, \inf_M u \geq m >0,  $$

for all $ u $ solution of $ (\ref{(1)}) $.
\end{Corollary}

\section{Proof of the results}
\underbar {\it Proof of theorem 1.1:}
\bigskip

Let $ x_0 $ be a point of $ M $. We want to prove a uniform estimate around $ x_0 $.

\bigskip

Let $ (u_i)_i $ be a sequence of solutions to:

$$ \Delta u_i + h u_i =8{u_i}^3, \,\, u_i>0, $$

We argue by contradiction, we assume that the $ \sup \times \inf $ is not bounded.

\bigskip

$ \forall \,\, c,R >0 \,\, \exists \,\, u_{c,R} $ solution to $ (1) $ such that:

\be R^2 (\sup_{B(x_0,R)} u_{c,R})^{1/3}\times \inf_{B(x_0,2R)} u_{c,R} \geq c, \label{(4)}\ee

\begin{Proposition}{\it (blow-up analysis)} 

\smallskip

There is a sequence of points $ (y_i)_i $, $ y_i \to x_0 $  and two sequences of positive real numbers $ (l_i)_i, (L_i)_i $, $ l_i \to 0 $, $ L_i \to +\infty $, such that if we set $ v_i(y)=\dfrac{u_i[\exp_{y_i}(y/[u_i(y_i)])]}{u_i(y_i)} $, we have:

$$ 0 < v_i(y) \leq  \beta_i \leq 2, \,\, \beta_i \to 1. $$

$$  v_i(y)  \to \dfrac{1}{1+{|y|^2}}, \,\, {\rm uniformly \,\, on \,\, compact \,\, sets \,\, of } \,\, {\mathbb R}^4 . $$

$$ l_i^2 (u_i(y_i))^{1/3} \min_M u_i  \to +\infty. $$

\end{Proposition}

\underbar {\bf Proof:}

\bigskip

We use the hypothesis $ (\ref{(4)}) $, we take two sequences, $ R_i>0, R_i \to 0 $ and $ c_i \to +\infty $, such that,

\be {R_i}^2 (\sup_{B(x_0,R_i)} u_i)^{1/3}  \times \inf_{B(x_0,2R_i)} u_i\geq c_i \to +\infty, \label{(5)}\ee

Let, $ x_i \in  { B(x_0,R_i)} $, such that $ \sup_{B(x_0,R_i)} u_i=u_i(x_i) $ and $ s_i(x)=[R_i-d(x,x_i)] u_i(x), x\in B(x_i, R_i) $. Then, $ x_i \to x_0 $.

\bigskip

We have:

$$ \max_{B(x_i,R_i)} s_i(x)=s_i(y_i) \geq s_i(x_i)={R_i} u_i(x_i)\geq \sqrt {c_i}  \to + \infty. $$ 

We set :

$$ l_i=R_i-d(y_i,x_i),\,\, \bar u_i(y)= u_i [\exp_{y_i}(y)],\,\,  v_i(z)=\dfrac{u_i [ \exp_{y_i}\left ( z/[u_i(y_i)] \right )] } {u_i(y_i)}. $$

Clearly, we have, $ y_i \to x_0 $. We obtain:

$$ L_i= \dfrac{l_i}{(c_i)^{1/4}} [u_i(y_i)]=\dfrac{[s_i(y_i)]}{c_i^{1/4}}\geq \dfrac{c_i^{1/2}}{c_i^{1/4}}=c_i^{1/4}\to +\infty. $$

\bigskip

If $ |z|\leq L_i $, then $ y=\exp_{y_i}[z/ [u_i(y_i)]] \in B(y_i,\delta_i l_i) $ with $ \delta_i=\dfrac{1}{(c_i)^{1/4}} $ and $ d(y,y_i) < R_i-d(y_i,x_i) $, thus, $ d(y, x_i) < R_i $ and, $ s_i(y)\leq s_i(y_i) $. We can write,

$$ u_i(y) [R_i-d(y,y_i)] \leq u_i(y_i) l_i. $$

But, $ d(y,y_i) \leq \delta_i l_i $, $ R_i >l_i$ and $ R_i-d(y, y_i) \geq R_i-\delta_i l_i>l_i-\delta_i l_i=l_i(1-\delta_i) $, hence, we obtain,

$$ 0 < v_i(z)=\dfrac{u_i(y)}{u_i(y_i)} \leq \dfrac{l_i}{l_i(1-\delta_i)} \leq 2. $$

We set, $ \beta_i= \dfrac{1}{1-\delta_i} $, clearly $ \beta_i \to 1 $.

\bigskip

The function $ v_i $ satisfies the following equation:

\be -g^{jk}(z)\partial_{jk} v_i-\partial_k \left [ g^{jk}\sqrt { |g| } \right ](z)\partial_j v_i+ \dfrac{h(z)}{[u_i(y_i)]^2} v_i= 8{v_i}^3 \label{(6)}\ee
\bigskip

We use Ascoli and Ladyzenskaya theorems to obtain the local uniform convergence (on every compact set of $ {\mathbb R}^4 $) of $ ( v_i)_i $ to $ v $ solution on $ {\mathbb R}^4 $ to: 

$$ \Delta v = 8v^3, \,\, v(0)=1,\,\, 0 \leq v\leq 1\leq 2, $$

By the maximum principle, we have $ v>0 $ on $ {\mathbb R}^n $. According to Caffarelli-Gidas-Spruck result (see [10]), we have, $ v(y)= \dfrac{1}{1+{|y|^2}} $.

\bigskip

\underbar {\it Polar Geodesic Coordinates}

\bigskip

Let $ u $ be a function on $ M $. We set $ \bar u(r,\theta)=u[\exp_x(r\theta)] $. We denote $ g_{x,ij} $ the local expression of the metric $ g $ in the exponential chart centered at $ x $.

\bigskip

We set,

$$ w_i(t,\theta)=e^t\bar u_i(e^t,\theta) = e^{t}u_i[\exp_{y_i}(e^t\theta)] , $$

$$ a(y_i,t,\theta)=\log J(y_i,e^t,\theta)=\log [\sqrt { det(g_{y_i, ij})}]. $$ 

We can write the Laplace-Beltrami operator in polar geodesic coordinates:

\be -\Delta u =\partial_{rr} \bar u+\dfrac{3}{r} \partial_r \bar u+\partial_r [\log  J(x,r,\theta)] \partial_r \bar u-\dfrac{1}{r^2}\Delta_{\theta } \bar u . \label{(7)}\ee

We deduce the two following lemmas:

\bigskip

\begin{Lemma}

\smallskip

The function $ w_i $ is a solution to:

\be  -\partial_{tt} w_i-\partial_t a \partial_t w_i-\Delta_{\theta }w_i+c w_i=8w_i^3, \label{(8)}\ee

with

 $$ c = c(y_i,t,\theta)= 1+\partial_t a + h e^{2t}, $$ 

\end{Lemma}

\underbar {\bf Proof:}

\smallskip

We write:

$$ \partial_t w_i=e^{2t}\partial_r \bar u_i+ w_i,\,\, \partial_{tt} w_i=e^{3t} \left [\partial_{rr} \bar u_i+\dfrac{3}{e^t}\partial_r \bar u_i \right ]+ w_i. $$

$$ \partial_t a =e^t\partial_r \log J(y_i,e^t,\theta), \partial_t a \partial_t w_i=e^{3t}\left [ \partial_r \log J\partial_r \bar u_i \right ]+ \partial_t a w_i.$$

Lemma 1 follows.

\bigskip

Let $ b_1(y_i,t,\theta)=J(y_i,e^t,\theta)>0 $. We can write:

$$ -\dfrac{1}{\sqrt {b_1}}\partial_{tt} (\sqrt { b_1} w_i)-\Delta_{ \theta }w_i+[c(t)+ b_1^{-1/2} b_2(t,\theta)]w_i=8{w_i}^3, $$

where, $ b_2(t,\theta)=\partial_{tt} (\sqrt {b_1})=\dfrac{1}{2 \sqrt { b_1}}\partial_{tt}b_1-\dfrac{1}{4(b_1)^{3/2}}(\partial_t b_1)^2 .$

\bigskip

We set,

$$ \tilde w_i=\sqrt {b_1} w_i. $$

\begin{Lemma}

\smallskip

The function $ \tilde w_i $ is a solution to:

$$ -\partial_{tt} \tilde w_i+\Delta_{ \theta } (\tilde w_i)+2\nabla_{\theta}(\tilde  w_i) .\nabla_{\theta} \log (\sqrt {b_1})+(c+b_1^{-1/2} b_2-c_2) \tilde w_i= $$

\be = 8\left (\dfrac{1}{b_1} \right ){\tilde w_i}^3, \label{(9)}\ee

where, $ c_2 $ is a function to be determined.

\end{Lemma}

\bigskip

\underbar {\bf Proof:}

We have: 

$$ -\partial_{tt} \tilde w_i-\sqrt {b_1} \Delta_{ \theta } w_i+(c+b_2) \tilde w_i= 8 \left (\dfrac{1}{b_1} \right ) {\tilde w_i}^3, $$

But,

$$ \Delta_{ \theta } (\sqrt {b_1} w_i)=\sqrt {b_1} \Delta_{ \theta } w_i-2 \nabla_{\theta} w_i .\nabla_{\theta} \sqrt {b_1}+ w_i \Delta_{ \theta }(\sqrt {b_1}), $$

and,

$$ \nabla_{\theta} (\sqrt {b_1} w_i)=w_i \nabla_{\theta} \sqrt {b_1}+ \sqrt {b_1} \nabla_{\theta} w_i, $$

we can write,

$$ \nabla_{\theta} w_i. \nabla_{\theta} \sqrt {b_1}=\nabla_{\theta}(\tilde  w_i) .\nabla_{\theta} \log (\sqrt {b_1})-\tilde w_i|\nabla_{\theta} \log (\sqrt {b_1})|^2, $$

we deduce,

$$  \sqrt {b_1} \Delta_{\theta } w_i= \Delta_{ \theta } (\tilde w_i)+2\nabla_{\theta}(\tilde  w_i) .\nabla_{\theta} \log (\sqrt {b_1})-c_2 \tilde w_i, $$

with $ c_2=[\dfrac{1}{\sqrt {b_1}} \Delta_{ \theta }(\sqrt{b_1}) + |\nabla_{\theta} \log (\sqrt {b_1})|^2] $. Lemma 2 is proved.

\bigskip

\underbar {\it The moving-Plane method:}

\bigskip

Let $ \xi_i $  be a real number,  we assume $ \xi_i \leq t $. We set $ t^{\xi_i}=2\xi_i-t $ and $ \tilde w_i^{\xi_i}(t,\theta)=\tilde w_i(t^{\xi_i},\theta) $. Set, $ \lambda_i=-\log u_i(y_i) $

\bigskip

\begin{Proposition}

\smallskip

We claim:

\be\tilde w_i(\lambda_i,\theta)-\tilde w_i(\lambda_i+4,\theta) \geq \tilde k>0, \,\, \forall \,\, \theta \in {\mathbb S}_{3}. \label{(10)}\ee

For all $ \beta >0 $, there exists $ c_{\beta} >0 $ such that:

\be \dfrac{1}{c_{\beta}} e^t \leq \tilde w_i(\lambda_i+t,\theta) \leq c_{\beta}e^t, \,\, \forall \,\, t\leq \beta, \,\, \forall \,\, \theta \in {\mathbb S}_{3}. \label{(11)}\ee

\end{Proposition}

\underbar {\bf Proof:}

\smallskip

As in [2], we have, $ w_i(\lambda_i, \theta)-w_i(\lambda_i+4,\theta) \geq k>0 $ for $ i $ large, $ \forall \,\, \theta $. We can remark that $ b_1(y_i,\lambda_i,\theta) \to 1 $ and $ b_1(y_i,\lambda_i+4,\theta) \to 1 $ uniformly in $ \theta $, we obtain the first claim of proposition 2.4. For the second claim we use proposition 2.1, see also [2].

\bigskip

We set:

\be \bar Z_i=-\partial_{tt} (...)+\Delta_{ \theta } (...)+2\nabla_{\theta}(...) .\nabla_{\theta} \log (\sqrt {b_1})+(c+b_1^{-1/2} b_2-c_2)(...) \label{(12)}\ee

{\bf Remark :} In the operator $ \bar Z_i $, we can remark that:

$$ c+b_1^{-1/2}b_2-c_2 \geq k'>0,\,\, {\rm for }\,\, t<<0, $$

we can apply the maximum principle and the Hopf lemma.

\bigskip

\underbar {\bf Goal:}

\bigskip

Like in [2], we have an elliptic second order operator. Here it is $ \bar Z_i $, the goal is to use the "moving-plane" method to have a contradiction. For this, we must have:

\be \bar Z_i(\tilde w_i^{\xi_i}-\tilde w_i) \leq 0, \,\, {\rm if} \,\, \tilde w_i^{\xi_i}-\tilde w_i \leq 0. \label{(13)}\ee

We write, $ \Delta_{\theta}=\Delta_{g_{y_i, e^t, {}_{{\mathbb S}_{n-1}}}} $. We obtain:

$$ \bar Z_i(\tilde w_i^{\xi_i}-\tilde w_i)= (\Delta_{g_{y_i, e^{t^{\xi_i}}, {}_{{\mathbb S}_{3}}}}-\Delta_{g_{y_i, e^{t}, {}_{{\mathbb S}_{3}}}}) (\tilde w_i^{\xi_i})+ $$

$$ +2(\nabla_{\theta, e^{t^{\xi_i}}}-\nabla_{\theta, e^t})(w_i^{\xi_i}) .\nabla_{\theta, e^{t^{\xi_i}}} \log (\sqrt {b_1^{\xi_i}})+ 2\nabla_{\theta,e^t}(\tilde w_i^{\xi_i}).\nabla_{\theta, e^{t^{\xi_i}}}[\log (\sqrt {b_1^{\xi_i}})-\log \sqrt {b_1}]+ $$ 

$$ +2\nabla_{\theta,e^t} w_i^{\xi_i}.(\nabla_{\theta,e^{t^{\xi_i}}}-\nabla_{\theta,e^t})\log \sqrt {b_1}- [(c+b_1^{-1/2} b_2-c_2)^{\xi_i}-(c+b_1^{-1/2}b_2-c_2)]\tilde w_i^{\xi_i} + $$

\be + 8 \left ( \dfrac{1}{b_1^{\xi_i}} \right )({\tilde w_i}^{\xi_i})^3-8\left ( \dfrac{1}{b_1} \right ) {\tilde w_i}^3. \label{(14)}\ee

Clearly, we have the following lemma:

\bigskip

\begin{Lemma}

$$ b_1(y_i,t,\theta)=1-\dfrac{1}{3} Ricci_{y_i}(\theta,\theta)e^{2t}+\ldots, $$

$$ R_g(e^t\theta)=R_g(y_i) + <\nabla R_g(y_i)|\theta > e^t+\dots . $$

\end{Lemma}

According to proposition 1 and lemma 3,

\bigskip

\begin{Proposition}

$$ \bar Z_i(\tilde w_i^{\xi_i}-\tilde w_i) \leq 8 {(b_1^{\xi_i})}[(\tilde w_i^{\xi_i})^3- \tilde w_i^3] + +C|e^{2t}-e^{2t^{\xi_i}}|(|\nabla_{\theta} {\tilde w_i}^{\xi_i}| + |\nabla_{\theta}^2(\tilde w_i^{\xi_i})|) + $$

\be +C|e^{2t}-e^{2t^{\xi_i}}| (|Ricci_{y_i}|+|h|) \tilde w_i^{\xi_i}+ C' w_i^{\xi_i} |e^{3t^{\xi_i}}-e^{3t}|. \label{(15)}\ee

\end{Proposition}

\underbar {\bf Proof of proposition 2.6:}

\bigskip

In polar geodesic coordinates (and the Gauss lemma):

\be g = dt^2+r^2{\tilde g}_{ij}^kd\theta^id\theta^j \,\, {\rm et } \,\, \sqrt { |{\tilde g}^k|}=\alpha^k(\theta) \sqrt {[det(g_{x,ij})]}, \label{(16)}\ee

where $ \alpha^k $ is the volume element of the unit sphere associated to $ U^k $.

\smallskip

We can write (with lemma 2.3):

$$ |\partial_t b_1(t)|+|\partial_{tt} b_1(t)|+|\partial_{tt} a(t)|\leq C e^{2t}, $$

and,

$$ |\partial_{\theta_j} b_1|+|\partial_{\theta_j,\theta_k} b_1|+\partial_{t,\theta_j}b_1|+|\partial_{t,\theta_j,\theta_k} b_1|\leq C e^{2t}, $$

But,

$$ \Delta_{\theta}=\Delta_{g_{y_i, e^t, {}_{{\mathbb S}_{3}}}}=-\dfrac{\partial_{\theta^l}[{\tilde g}^{\theta^l \theta^j}(e^t,\theta)\sqrt { |{\tilde g}^k(e^t,\theta)|}\partial_{\theta^j}]}{\sqrt {|{\tilde g}^k(e^t,\theta)|}} . $$

Then,

\be A_i:=\left [{ \left [ \dfrac{\partial_{\theta^l}({\tilde g}^{\theta^l \theta^j}\sqrt { |{\tilde g}^k|}\partial_{\theta^j})}{\sqrt {|{\tilde g}^k|}} \right ] }^{\xi_i}-\left [ \dfrac{\partial_{\theta^l}({\tilde g}^{\theta^l \theta^j}\sqrt { |{\tilde g}^k|}\partial_{\theta^j})}{\sqrt {|{\tilde g}^k|}} \right ] \right ](\tilde w_i^{\xi_i}) = B_i+D_i \label{(17)}\ee

where,

\be B_i=\left [ {\tilde g}^{\theta^l \theta^j}(e^{t^{\xi_i}}, \theta)-{\tilde g}^{\theta^l \theta^j}(e^t,\theta) \right ] \partial_{\theta^l \theta^j}\tilde w_i^{\xi_i}, \label{(18)}\ee

and,

\be D_i=\left [ \dfrac{\partial_{\theta^l}[{\tilde g }^{\theta^l \theta^j}(e^{t^{\xi_i}},\theta)\sqrt {| {\tilde g}^k|}(e^{t^{\xi_i}},\theta)  ]}{ \sqrt {| {\tilde g}^k|}(e^{t^{\xi_i}},\theta )  }            -\dfrac{\partial_{\theta^l} [{\tilde g }^{\theta^l \theta^j}(e^t,\theta)\sqrt {| {\tilde g}^k|}(e^t,\theta) ]}{ \sqrt {| {\tilde g}^k|}(e^t,\theta) } \right ] \partial_{\theta^j} \tilde w_i^{\xi_i}. \label{(19)}\ee

Clearly, we can choose $ \epsilon_1 >0 $ such that:

\be |\partial_r{\tilde g}_{ij}^k(x,r,\theta)|+|\partial_r\partial_{\theta^m}{\tilde g}_{ij}^k(x,r, \theta)| \leq C r,\,\, x\in B(x_0,\epsilon_1) \,\, r\in [0,\epsilon_1], \,\,\theta \in U^k.\label{(20)}\ee

finally,

\be A_i \leq C_k|e^{2t}-e^{2t^{\xi_i}}|\left [ |\nabla_{\theta} \tilde w_i^{\xi_i}| + |\nabla_{\theta}^2(\tilde w_i^{\xi_i})| \right ], \label{(21)}\ee

We take, $ C=\max \{ C_i, 1 \leq i\leq q \} $ and we use $ (\ref{(14)}) $. Proposition 2.6 is proved.

\bigskip

We have,

\be c(y_i,t,\theta)= 1 + \partial_t a + h e^{2t}, \label{(22)}\ee

\be b_2(t,\theta)=\partial_{tt} (\sqrt {b_1})=\dfrac{1}{2 \sqrt { b_1}}\partial_{tt}b_1-\dfrac{1}{4(b_1)^{3/2}}(\partial_t b_1)^2 , \label{(23)}\ee 

\be c_2=[\dfrac{1}{\sqrt {b_1}} \Delta_{ \theta }(\sqrt{b_1}) + |\nabla_{\theta} \log (\sqrt {b_1})|^2], \label{(24)}\ee

\bigskip

We assume that $ \lambda \leq \lambda_i+2=-\log u_i(y_i)+2 $, which will be choosen later.

\smallskip

We work on $ [\lambda,t_i] \times {\mathbb S}_3 $ with $ t_i = \log { l_i } \to -\infty $, $ l_i $ as in the proposition 1. For $ i $ large $ \log {l_i} >> \lambda_i + 2 $.

\bigskip

The functions  $ v_i $ tend to a radially symmetric function, then, $ \partial_{\theta_j} w_i^{\lambda} \to 0 $ if $ i \to +\infty $ and,

$$ \dfrac{\partial_{\theta_j}w_i^{\lambda }(t,\theta)}{w_i^{\lambda }}=\dfrac{e^{(n-2)[(\lambda -\lambda_i)+(\xi_i-t)]/2} e^{[(\lambda -\lambda_i)+(\xi_i-t)]}(\partial_{\theta_j} v_i)(e^{[(\lambda -\lambda_i)+(\lambda -t)]}\theta)}{e^{(n-2)[(\lambda-\lambda_i)+(\lambda-t)]/2}v_i[e^{(\lambda-\lambda_i)+(\lambda - t)}\theta]} \leq {\bar C_i}, $$

where $ \bar C_i $ does not depend on $ \lambda $ and tends to 0. We have also,

\be |\partial_{\theta} w_i^{\lambda }(t,\theta)|+|\partial_{\theta,\theta} w_i^{\lambda }(t,\theta)|\leq {\tilde C_i} w_i^{\lambda }(t,\theta), \,\, {\tilde C_i} \to 0. \label{(25)}\ee

and,

\be |\partial_{\theta} \tilde w_i^{\lambda }(t,\theta)|+|\partial_{\theta,\theta}  \tilde w_i^{\lambda }(t,\theta)|\leq {\tilde C_i} \tilde w_i^{\lambda }(t,\theta), \,\, {\tilde C_i} \to 0. \label{(26)}\ee

$ \tilde C_i $ does not depend on $ \lambda $.

Now, we set:

\be \bar w_i=\tilde w_i-\dfrac{\tilde m_i}{2} e^{2t}, \label{(27)}\ee

with, $ m_i=\dfrac{1}{2} u_i(x_i)^{1/3} \min_M u_i $. As in [2], we have,

\bigskip

\begin{Lemma}

\bigskip

There is $ \nu <0 $ such that for $ \lambda \leq \nu $ :

\be \bar w_i^{\lambda}(t,\theta)-\bar w_i(t,\theta) \leq 0, \,\, \forall \,\, (t,\theta) \in [\lambda,t_i] \times {\mathbb S}_{3}. \label{(28)}\ee

\end{Lemma}

Let $ \xi_i $ be the following real number,

$$ \xi_i=\sup \{ \lambda \leq \lambda_i+2, \bar w_i^{\xi_i}(t,\theta)-\bar w_i(t,\theta) \leq 0, \,\, \forall \,\, (t,\theta)\in [\xi_i,t_i]\times {\mathbb S}_{3} \}. $$

Like in [2], we use the previous lemma to show:

$$ \bar w_i^{\xi_i}-\bar w_i \leq 0 \Rightarrow \bar Z_i(\bar w_i^{\xi_i}-\bar w_i) \leq 0. $$

We have,

$$ \bar Z_i(\tilde w_i^{\xi_i}-\tilde w_i) \leq  8b_1^{\xi_i}[( \tilde w_i^{\xi_i})^3- \tilde w_i^3]+ O(1)(e^{2t}-e^{2t^{\xi_i}})+O(1) {\tilde w_i}^{\xi_i}(e^{2t}-e^{2t^{\xi_i}}). $$



$$ -\bar Z_i(e^{2t^{\xi_i}}-e^{2t})=(4-1-\partial_t a-he^{2t}+b_1^{-1/2}b_2-c_2)(e^{2t^{\xi_i}}-e^{2t}) \leq c_3(e^{2t^{\xi_i}}-e^{2t}) $$

Thus,

$$ \bar Z_i(\bar w_i^{\xi_i}-\bar w_i) \leq 8 b_1^{\xi_i}[( \tilde w_i^{\xi_i})^3-{\tilde w_i}^3]+ (c_3m_i-c_4) (e^{2t^{\xi_i}}-e^{2t} ). $$

with, $ c_3, c_4 >0 $.

But,

$$  0 < {\tilde w_i}^{\xi_i} \leq 2e, \,\,\, \tilde w_i \geq \dfrac{m_i}{2} e^{2t} \,\,\, {\rm and} \,\,\, {\tilde w_i}^{\xi_i}-{\tilde w_i} \leq \dfrac{m_i}{2} ( e^{2t^{\xi_i}}-e^{2t}), $$

and,

\be ({\tilde w_i}^{\xi_i})^3-{\tilde w_i}^3=({\tilde w_i}^{\xi_i}-\tilde w_i)[ ({\tilde w_i}^{\xi_i})^2+{\tilde w_i}^{\xi_i} {\tilde w_i}+ {\tilde w_i}^2] \leq ({\tilde w_i}^{\xi_i}-{\tilde w_i}) ({\tilde w_i^{\xi_i}})^2 + ({\tilde w_i}^{\xi_i}-{\tilde w_i}) \dfrac{m^2 e^{2t}}{4} + ({\tilde w_i}^{\xi_i}-{\tilde w_i}) \dfrac{m}{2} e^t {\tilde w_i}^{\xi_i}, \label{(29)}\ee

then,

\be \bar Z_i(\bar w_i^{\xi_i}-\bar w_i) \leq   [ [\dfrac{a m_i^3}{16}-O(1)] + [\dfrac{a m_i^2}{8}- O(1)] e^t {\tilde w_i}^{\xi_i}] ( e^{2t^{\xi_i}}- e^{2t}) \leq 0. \label{(30)}\ee

If we use the maximum principle and the Hopf lemma, we obtain (as in [2]):

$$ \max_{\theta \in {\mathbb S}_3} w_i(t_i,\theta) \leq \min_{\theta \in {\mathbb S}_3} w_i(2\xi_i-t_i, \theta), $$

we can write (using  proposition 2):

\be  l_i [u_i(y_i)]^{1/3} \min_M u_i \leq c, \label{(31)}\ee

\section{Metric and Laplacian}

In this section, we give some remarks on Polar Geodesic coordinates and the Laplacian in these coordinates. First by using the Jacobi Fields we can have an expansion of the metric in geodesic coordinates, we can extend this result to polar geodesic coordinates.

{\underbar {\bf Estimate of the metric in Polar Coordinates.}}

\bigskip

Let us consider a riemannian manifold (not necessarily compact) $ (M,g) $. We set $ g_{x,ij} $ the components of the metric in the exponential chart centered at $ x $.

By the Gauss formula we have:

$$ g=ds^2=dt^2+g_{ij}^k(r,\theta)d\theta^id\theta^j=dt^2+r^2{\tilde g}_{ij}^k(r,\theta)d\theta^id\theta^j=g_{x,ij}dx^idx^j, $$

in polar chart centered at $ x $, $ ]0,\epsilon_0[\times U^k $, with $ ( U^k, \psi) $ a chart of the unit sphere $ {\mathbb S}_{n-1} $.

We can write the volume lement as:

$$ dV_g=r^{n-1}\sqrt {|{\tilde g}^k|}dr d \theta^1 \ldots d \theta^{n-1} = \sqrt {[det(g_{x,ij})]}dx^1 \ldots dx^n, $$

thus,

$$ dV_g=r^{n-1} \sqrt { [det(g_{x,ij})]}[\exp_x(r\theta)]\alpha^k(\theta)dr d\theta^1 \ldots d \theta^{n-1} , $$

with, $ \alpha^k $ is such that $ d\sigma_{{\mathbb S}_{n-1}}=\alpha^k(\theta) d \theta^1 \ldots d \theta^{n-1} . $ (volume element of the sphere in the chart $ (U^k,\psi) $ ).

\bigskip

Thus,

$$ \sqrt { |{\tilde g}^k|}=\alpha^k(\theta) \sqrt {[det(g_{x,ij})]}, $$

\underbar {\bf Proposition:} Let us  consider $ x_0 \in M $,  there is $ \epsilon_1>0 $  and $ U^k $,  such that:

$$ |\partial_r{\tilde g}_{ij}^k(x,r,\theta)|+|\partial_r\partial_{\theta^m}{\tilde g}_{ij}^k(x,r, \theta)| \leq C r,\,\, \forall \,\, x\in B(x_0,\epsilon_1) \,\, \forall \,\, r\in [0,\epsilon_1], \,\, \forall \,\, \theta \in U^k.$$

and,

$$ |\partial_r|{\tilde g}^k|(x,r,\theta)|+\partial_r \partial_{\theta^m} |{\tilde g}^k|(x,r,\theta)\leq C r,\,\, \forall \,\, x\in B(x_0,\epsilon_1) \,\, \forall \,\, r\in [0,\epsilon_1], \,\, \forall \,\, \theta \in U^k. $$

and,

$$ |\partial_r \partial_{\theta_m}\left [\dfrac{\sqrt {|{\tilde g}^k|} }{\alpha^k(\theta) } \right ] (x,r,\theta)|+|\partial_r \partial_{\theta_m}\partial_{\theta_{m'}} \left [\dfrac{\sqrt {|{\tilde g}^k|} }{\alpha^k(\theta) } \right ] (x,r,\theta)|\leq C r \,\,\forall \,\, (x,r,\theta)\in B(x_0,\epsilon_1)\times [0,\epsilon_1]\times U^k. $$

\underbar{\bf Proof:}

\bigskip

Next, we use Einstein convention:

\smallskip

First, we consider a chart $ (\Omega,\phi) $ in $ x_0 $, such that $ \bar \Omega $ is compact, we can assume it normal at $ x_0 $.

\bigskip

According to lemma 2.3.7 of [He], for all $ x_0 \in M $ there is $ \epsilon_0>0 $ such that the application $ u:(x,v)\to \exp_x(v) $ on $ B(x_0,\epsilon_0) \times B(0,\epsilon_0) $ in $ M $ is $ {C}^{\infty} $  and for all $ x\in B(x_0,\epsilon) $ the application $ v\to \exp_x(v) $ is a diifeomorphism of $ B(0,\epsilon_0) $ in $ B(x_0,\epsilon_0) $ with $ \exp_x[\partial B(0,\mu)]=\partial B(x,\mu) $, $ \mu \leq \epsilon_0 $.  Without loss of generality  we can assume that $ B(x,\epsilon_0) \subset \subset \Omega $ for all $ x \in B(x_0,\epsilon_0) $.

\bigskip

Thus, we have for all $ x\in B(x_0,\epsilon_0) $, $ [B(x,\epsilon_0), \exp_x^{-1}] $ is a normal chart in $ x $.(In this case we can define polar coordinates). For all $ x\in B(x_0,\epsilon_0) $:

$$ g_{x,ij}(z)=g(z)(\partial_{z^i,x},\partial_{z^j,x}), $$

with $ \partial_{z^i,x} $  is the canonical vector field in the exponential chart.

\smallskip

We set $ a_i^k(z,x)=\dfrac{(\phi o \exp_x)^k}{\partial z^i}[\exp_x^{-1}(z)] $, then $ \partial_{z^i,x}=a_i^k(z,x) \partial_{u^k,\phi}, $

\smallskip

with $ \partial_{u^k,\phi} $, the canonical vector field with respect to chart $ (\Omega, \phi) $,  this vector field do not depends on $ x $ and the functions  $ a_i^k $ are regular of $ z $ and $ x $. We obtain,

$$ g_{x,ij}(z)=g(z)[a_i^k(z,x) \partial_{u^k,\phi}; a_j^l(z,x)\partial_{u^l,\phi}]=a_i^k(z,x)a_j^l(z,x)g_{kl}(z), $$

with $ g_{kl} $ the component of $ g $  in the chart $ ( \Omega, \phi ) $.

\bigskip

We have $ z=\exp_x(y) $, $ y \in B(0,\epsilon_0) \subset {\mathbb R}^n $ and $ y=r\theta $ in polar coordinates, thus $ (x,r,\theta) \to g_{x,ij}[\exp_x(r\theta)] $ is $ { C}^{\infty} $ of $ x, r $ et $ \theta $.

\bigskip

We have, by definition, $ g_{ij}^k(r\theta)=g_{[\exp_x(r\theta)]}(\partial_{\theta^i,x},\partial_{\theta^j,x}) $, (canonical vector fields).

\bigskip

We can write,

$$ \partial_{\theta^i,x}=r b_i^k(\theta)\partial_{z^k,x} , $$

with $ b_i^j $ regular. (Note, that here, we can use the function $ \theta^i $, as regular function in the chart of the sphere. With this procedure we can deduce the component of the metric in polar coordinates ( a good expansion)). Thus,

$$ g_{ij}^k(r,\theta)=r^2g[\exp_x(r\theta)](b_i^k\partial_{z^k,x},b_j^l\partial_{z^l,x}), $$

Then,

$$ g_{ij}^k(r,\theta)=r^2b_i^k(\theta)b_j^l(\theta)g_{x,kl}[\exp_x(r\theta)]. $$

Thus, the functions $ {\tilde g}_{ij}^k : (x,r,\theta) \to b_j^kb_j^lg_{x,kl}[\exp_x(r\theta)] $ is regular of $ x,r $ and $ \theta $. We have,

$$ \partial_r {\tilde g}_{ij}^k(x,0,\theta)=b_i^k(\theta)b_j^k(\theta)c^m(\theta)\partial_mg_{x,kl}(x)=0 , $$

because the exponential chart is normal in $ x $ and $ g_{x,kl} $ are the component of $ g $ in this chart. We have:

$$ \partial_{\theta^m}{\tilde g}_{ij}^k(x,r,\theta)={\tilde b}_i^k(\theta) {\tilde b}_j^l(\theta)g_{x,kl}[\exp_x(r\theta)]+r{\bar b}_i^k(\theta){\bar b}_j^l(\theta){\bar c}_m^s(\theta)\partial_sg_{x,kl}[\exp_x(r\theta)], $$

We also have:

$$ \partial_r \partial_{\theta^m} {\tilde g}_{ij}^k(x,r,\theta)=u_{ijmr}^{klq}(\theta)\partial_q g_{x,kl}[\exp_x(r\theta)]+v_{ijmr}^{kl}(\theta)w^s(t)\partial_sg_{x,kl}[\exp_x(r\theta)]+ $$

$$ + rh_{ijrm}^{klst}(\theta)\partial_{st}g_{x,kl}[\exp_x(r\theta)]. $$

Thus,

$$ \partial_r \partial_{\theta^m} {\tilde g}_{ij}^k(x,0,\theta)=0, \,\, \forall \,_, x\in B(x_0,\epsilon_0), \,\, \forall \,\, \theta \in U^k . $$

Thus, we obtain:

$$ \partial_r {\tilde g}_{ij}^k(x,0,\theta)=\partial_r \partial_{\theta^m}{\tilde g}_{ij}^k(x,0,\theta)=0,\,\, \forall \,\, x\in B(x_0,\epsilon_0), \,\, \forall \,\, \theta \in U^k. \qquad \qquad (*) $$

\bigskip

Because $ \sqrt { |{\tilde g}^k|}=\alpha^k(\theta) \sqrt {[det(g_{x,ij})]} $, we deduce,

$$ \partial_r (\log \sqrt {|{\tilde g}^k|})=\partial_r [\log (\sqrt { [det(g_{x,ij})]})]. $$

We use the definition of the determinant

$$ det [g_{x,ij}][\exp_x(r\theta)]=\Sigma \,\, \Pi \,\, g_{x,kl}[\exp_x(r\theta)], $$

Thus,

$$ \partial_r det[g_{x,ij}](x)=\Sigma\,\, \Pi \,\, [g_{x,kl}(x)] a^s(\theta) \partial_sg_{x,mn}(x) = 0, $$

because the exponential chart is normal at $ x $.

\bigskip

Finaly

$$ \partial_r |{\tilde g}^k|(x,0,\theta)=0, \,\, \forall \,\, x\in B(x_0,\epsilon_0), \,\, \forall \,\, \theta \in U^k . $$

We also have $ \partial_r\partial_{\theta^m} {\tilde g}_{ij}^k(x,0,\theta) $, to prove that,

$$ \partial_r \partial_{\theta^m} |{\tilde g}^k|(x,0,\theta)=0. $$

If we set $ D_m=\partial_{\theta_m} \dfrac{\sqrt {|{\tilde g}^k| }}{\alpha^k(\theta) }(x,r,\theta) $,

$$ D_m=r \Sigma \, \Pi \,\, \beta_m^l(\theta)\partial_l g_{x, ij}[\exp_x(r\theta)] g_{x,ij}[\exp_x(r\theta)], $$

\bigskip

We have $ D_m(x,0,\theta)=0 $ and,

$$ \partial_r D_m(x,0,\theta)=\lim_{r\to 0} (D_m /r)(x,r,\theta)=0, $$

Thus

$$ \partial_r \partial_{\theta_m}\left [\dfrac{\sqrt {|{\tilde g}^k|} }{\alpha^k(\theta) } \right ] (x,0,\theta)=0, $$

$$ \partial_{\theta_{m'}}D_m= r\Sigma \, \Pi\,\, \partial_{m'}\beta_m^l \partial_l g_{x,ij} g_{x,ij}+r^2\Sigma \, \Pi \,\, \beta_m^l \beta_{m'}^{l'} \partial_{ll'}g_{x,ij}g_{x,ij}+r^2 \Sigma \, \Pi \,\, \beta_m^l \beta_{m'}^{l'}\partial_lg_{x,ij}\partial_{l'}g_{x,ij}g_{x,ij}, $$

but, $ \partial_{\theta_{m'}} D_m(x,0,\theta)=0 $, we have,

$$ \partial_r \partial_{\theta_{m'}} D_m(x,0,\theta)=\lim_{r\to 0} [\partial_{\theta_{m'}} D_m/r]=0, $$

Thus,

$$ \partial_r \partial_{\theta_m}\partial_{\theta{m'}} \left [\dfrac{\sqrt {|{\tilde g}^k|} }{\alpha^k(\theta) } \right ] (x,0,\theta)=0, $$

Finaly

$$  \partial_r {\tilde g}_{ij}^k(x,0,\theta)=\partial_r \partial_{\theta^m}{\tilde g}_{ij}^k(x,0,\theta)=0\,\, \forall \,\, x\in B(x_0,\epsilon_0), \,\, \forall \,\, \theta \in U^k. \qquad \qquad (**) $$

$$ \partial_r |{\tilde g}^k|(x,0,\theta)=\partial_r \partial_{\theta^m} |{\tilde g}^k|(x,0,\theta)=0\,\, \forall \,\, x\in B(x_0,\epsilon_0), \,\, \forall \,\, \theta \in U^k. \qquad \qquad (***) $$

$$ \partial_r \partial_{\theta_m}\left [\dfrac{\sqrt {|{\tilde g}^k|} }{\alpha^k(\theta) } \right ] (x,0,\theta)=\partial_r \partial_{\theta_m}\partial_{\theta{m'}} \left [\dfrac{\sqrt {|{\tilde g}^k|} }{\alpha^k(\theta) } \right ] (x,0,\theta)=0 \,\, \forall \,\, x\in B(x_0,\epsilon_0), \,\, \forall \,\, \theta \in U^k. \,\, (****) $$

We can reduce the open set $ U^k $ to have these estimates and use the uniform continuity to have the estimates.

\bigskip

{\underbar { \bf The Laplacian in polar coordinates}}

\bigskip

We can write in  $ [0,\epsilon_1]\times U^k $,

$$ -\Delta = \partial_{rr}+\dfrac{n-1}{r}\partial_r+ \partial_r [\log \sqrt { |{\tilde g^k|}] }\partial_r+\dfrac{1}{r^2 \sqrt {|{\tilde g}^k|}}\partial_{\theta^i}({\tilde g}^{\theta^i \theta^j}\sqrt { |{\tilde g}^k|}\partial_{\theta^j}) . $$

On a,

$$ -\Delta = \partial_{rr}+\dfrac{n-1}{r}\partial_r+ \partial_r \log J(x,r,\theta)\partial_r+ \dfrac{1}{r^2 \sqrt {|{\tilde g}^k|}}\partial_{\theta^i}({\tilde g}^{\theta^i \theta^j}\sqrt { |{\tilde g}^k|}\partial_{\theta^j}) . $$

We can write the Laplacian (radial and angular decomposition, see S. Lang book),

$$ -\Delta = \partial_{rr}+\dfrac{n-1}{r} \partial_r+\partial_r [\log J(x,r,\theta)] \partial_r-\Delta_{{\mathbb S}_r(x)}, $$

with $ \Delta_{{\mathbb S}_r(x)} $ the Laplacian on $ {\mathbb S}_r(x) $. 

\bigskip

Locally we can write $ \Delta_{{\mathbb S}_r(x)} $ as: 

$$ -\Delta_{{\mathbb S}_r(x)}=\dfrac{1}{r^2 \sqrt {|{\tilde g}^k|}}\partial_{\theta^i}({\tilde g}^{\theta^i \theta^j}\sqrt { |{\tilde g}^k|}\partial_{\theta^j}) . $$

\bigskip

We set: $ L_{\theta}(x,r)(...)=r^2\Delta_{{\mathbb S}_r(x)}(...)[\exp_x(r\theta)] $.

\bigskip

The operator $ L_{\theta}(x,r) $ act on $ { C}^2({\mathbb S}_{n-1}) $ functions globaly and not depends on the chart on $ {\mathbb S}_{n-1} $  and locally we can write it  as:

$$ L_{\theta}(x,r) =-\dfrac{1}{\sqrt {|{\tilde g}^k|}}\partial_{\theta^i}[{\tilde g}^{\theta^i \theta^j}\sqrt { |{\tilde g}^k|}\partial_{\theta^j}] . $$

We have 

$$ \Delta = \partial_{rr}+\dfrac{n-1}{r} \partial_r+\partial_r [ J(x,r,\theta)] \partial_r - \dfrac{1}{r^2} L_{\theta}(x,r) . $$

 We set, for $ u $  a function on $ M $,  $ \bar u(r,\theta)=uo\exp_x(r\theta) $ in polar coordnates centered in $ x $:

$$ -\Delta u=\partial_{rr} \bar u+\dfrac{n-1}{r} \partial_r \bar u+\partial_r [ J(x,r,\theta)] \partial_r \bar u-\Delta_{{\mathbb S}_r(x)}(u_{|{\mathbb S}_r(x)})[\exp_x(r\theta)], $$

$$ r^2\Delta_{{\mathbb S}_r(x)}(u_{|{\mathbb S}_r(x)})[\exp_x(r\theta)]=-\dfrac{1}{ \sqrt {|{\tilde g}^k|}}\partial_{\theta^i} \left [ {\tilde g}^{\theta^i \theta^j}\sqrt { |{\tilde g}^k|}\partial_{\theta^j} [uo\exp_x(r\theta)] \right ] =L_{\theta}(x,r)\bar u . $$

Thus,

$$ -\Delta u =\partial_{rr} \bar u+\dfrac{n-1}{r} \partial_r \bar u+\partial_r [ J(x,r,\theta)] \partial_r \bar u-\dfrac{1}{r^2}L_{\theta}(x,r)\bar u . $$

\underbar {\bf Lemma:}

\smallskip

The operator $ L_{\theta}(x,r) $ is a Laplacian on $ {\mathbb S}_{n-1} $ for a particular metric which depends on $ r $.

\bigskip

\underbar {\bf Proof on the Lemma:}

\bigskip

We have $ \Delta_{{\mathbb S}_r(x)}=\Delta_{i_{x,r}^{*}(g)} $ with $ i_{x,r} $ the identity map from $ {\mathbb S}_{r}(x) $ in $ M $ and $ i_{x,r}^*(g)=\tilde g $ the induced metric on the submanifold $ {\mathbb S}_r(x) $ of $ M $.

\bigskip

The map $ \exp_x $ induce a diffeomorphism from $ {\mathbb S}_r(x) $ into $ {\mathbb S}_{n-1}^r $, we have:

$$ \Delta_{{\mathbb S}_r(x)} u_{| {\mathbb S}_r(x)}=\Delta_{i_{x,r}^*(g)} u_{| {\mathbb S}_r(x)}=\Delta_{\tilde g} u_{| {\mathbb S}_r(x)}=\Delta_{\exp_x^*(\tilde g), {\mathbb S}_{n-1}^r} uo\exp_x(v), $$

Let us consider $ \tilde z $ the map from $ {\mathbb S}_{n-1} $ into $ {\mathbb S}_{n-1}^r $  defined by $ \theta \to r \theta $. Thus,

 $$ \Delta_{\exp_x^*(\tilde g)} uo\exp_x(v) =\Delta _{\tilde z^* [\exp_x^*(\tilde g)]} uo\exp_x(r\theta), $$

For a chart $ (\psi, U^k) $ on $ {\mathbb S}_{n-1} $, the polar chart in $ x $ is, $ (\phi_0, [0,\epsilon_0]\times U^k) $, where $ \phi_0= \exp_xo\tilde zo\psi^{-1} $. The $ g_{jl}^k $ are, by definition:

$$ g_{jl}^k(r,\theta)=g_{\exp_x(r\theta)}(\partial_{j,\phi_0}, \partial_{l,\phi_0}), $$

with, $ \partial_{j,\phi_0}, \partial_{l,\phi_0} $ the canonical vector fields for the chart $ (\phi_0, [0,\epsilon_0]\times U^k ) $.

\smallskip

By definition, take $ h $  a function defined in a neighborhood of $ x $:

$$ \partial_{j,\phi_0}(h)=\dfrac{\partial (ho\phi_0)}{\partial \bar \theta_j}=\dfrac{\partial (ho\exp_xo\tilde zo\psi^{-1})}{\partial \bar \theta_j}, $$ 

with $ \psi(\theta)=(\bar \theta_1,\ldots,\bar \theta_{n-1}) $ and $ \theta_0 = r $ (angular and radial derivations).

We look to the angular components, $ \partial_{j,\phi_0} $ with $ 1 \leq j\leq n-1 $, $ j \in {\mathbb N} $,thus,

$$ \partial_{j,\phi_0}(h)=\dfrac{\partial (hoi_{x,r}o\exp_xo\tilde zo\psi^{-1})}{\partial \bar \theta_j}, $$

$$ \dfrac{\partial (ho\phi_0)}{\partial \bar \theta_j}=\dfrac{\partial (hoi_{x,r}o\phi_0)}{\partial \bar \theta_j}, $$

Thus, if write $ \bar \partial_j, \bar \partial_l $ as a canonical vector fields of the unit sphre in the chart $ (\psi, U^k) $,  we obtain,  we use $ d $ for the differential,

$$ \partial_{j,\phi_0}=d(i_{x,r}o\exp_xo\tilde z)(\bar \partial_j), $$

Thus,

$$ g_{jl}^k(r,\theta)=g_{\exp_x(r\theta)}[d(i_{x,r}o\exp_xo\tilde z(\bar \partial_j), d(i_{x,r}o\exp_xo\tilde z(\bar \partial_l)], $$

we use the definition of the pull-back:

$$ g_{jl}^k(r,\theta)=\tilde z^*[\exp_x^*[i_{x,r}^*(g)]]_{\theta}(\bar \partial_j,\bar \partial_l)=\tilde z^*[\exp_x^*(\tilde g)](\bar \partial_j,\bar \partial_l), $$

Thus  we have the component $ j,l $ of the mteric $ \tilde z^* [\exp_x^*(\tilde g)] $ on the unit sphere $ {\mathbb S}_{n-1} $.

Finaly we have locally and globaly:

$$ -\Delta_{\tilde z^*[\exp_x^*(\tilde g)], {\mathbb S}_{n-1}} uo\exp_x(r\theta)=\dfrac{1}{r^2 \sqrt {|{\tilde g}^k|}}\partial_{\theta^i}[{\tilde g}^{\theta^i \theta^j}\sqrt { |{\tilde g}^k|}\partial_{\theta^j} uo\exp_x(r\psi^{-1})], $$

Now, if we consider the metric on $ {\mathbb S}_{n-1} $ defined by, 

$$ g_{{}_{x,r,{\mathbb S}_{n-1}}}=r^{-2} \tilde z^*[\exp_x^*(\tilde g)], $$

this metric is well-defined and is such that:

$$ {g_{{}_{x,r,{\mathbb S}_{n-1}}}}_{jl}=r^{-2}g_{jl}^k=\tilde g_{jl}^k, $$

\underbar{Bibliography for this section:}

\bigskip

[Au] T. Aubin. Some Nonlinear Problems in Riemannian Geometry. Springer-Verlag 1998.

[He] E. Hebey, Analyse non lineaire sur les Vari\'et\'es, Editions Diderot.

[L] S. Lang. Differentiable and Riemannian manifolds.

\end{document}